\input amstex\documentstyle{amsppt}  
\pagewidth{12.5cm}\pageheight{19cm}\magnification\magstep1
\topmatter
\title Generic character sheaves on disconnected groups and character values\endtitle
\author G. Lusztig\endauthor
\address{Department of Mathematics, M.I.T., Cambridge, MA 02139}\endaddress
\thanks{Supported in part by the National Science Foundation.}\endthanks
\endtopmatter   
\document

\define\dZ{\dot Z}

\define\uw{\un w}
\define\uD{\un D}

\define\si{\sim}
\define\wt{\widetilde}

\define\qua{\quad}

\define\tcl{\ti\cl}

\define\lb{\linebreak}

\define\part{\partial}
\define\em{\emptyset}

\define\m{\mapsto}
\define\do{\dots}

\define\lra{\leftrightarrow}

\define\T{\times}
\define\ti{\tilde}
\define\nl{\newline}
\redefine\i{^{-1}}

\define\un{\underline}

\define\ot{\otimes}
\define\bbq{\bar{\QQ}_l}

\define\Ad{\text{\rm Ad}}

\define\tr{\text{\rm tr}}

\define\a{\alpha}

\redefine\c{\chi}

\define\e{\epsilon}

\define\ph{\phi}
\define\ps{\psi}

\redefine\t{\tau}
\define\th{\theta}

\redefine\l{\lambda}

\define\x{\xi}

\redefine\G{\Gamma}

\define\kk{\bold k}

\define\uu{\bold u}

\define\ww{\bold w}

\define\FF{\bold F}

\define\NN{\bold N}

\define\QQ{\bold Q}

\define\ZZ{\bold Z}

\define\cb{\Cal B}

\define\ce{\Cal E}

\define\ch{\Cal H}

\define\cl{\Cal L}

\define\cn{\Cal N}
\define\co{\Cal O}

\define\tx{\ti x}
\define\ty{\ti y}

\define\tF{\ti F}

\define\tT{\ti T}

\define\che{\check}
\define\cha{\che{\a}}

\define\DL{DL}
\define\DM{DM}
\define\GF{L1}
\define\RC{L2}
\define\CDG{L3}
\define\MA{M}
\define\SH{S}

\head Introduction\endhead
The theory of character sheaves \cite{\CDG} on a reductive group $G$ over an algebraically
closed field and the theory of irreducible characters of $G$ over a finite field are two 
parallel theories; the first one is geometric (involving intersection cohomology complexes
on $G$), the second one involves functions on the group of rational points of $G$. In the 
case where $G$ is connected, a bridge between the two theories was constructed in 
\cite{\GF} and strengthened in \cite{\RC}, \cite{\SH}. In this paper we begin the 
construction of the analogous bridge in the general case, extending the method of 
\cite{\GF}. Here we restrict ourselves to character sheaves which are "generic" (in 
particular their support is a full connected component of $G$) and show how such character
sheaves are related to characters of representations (see Theorem 1.2).

\head Contents\endhead
1. Statement of the Theorem.

2. Constructing representations of $G^F$.

3. Proof of Theorem 1.2.

\head 1. Statement of the Theorem\endhead
\subhead 1.1\endsubhead
Let $\kk$ be an algebraic closure of a finite field $\FF_q$. Let $G$ be a reductive 
algebraic group over $\kk$ with identity component $G^0$ such that $G/G^0$ is cyclic,
generated by a fixed connected component $D$. We assume that $G$ has a fixed 
$\FF_q$-rational structure with Frobenius map $F:G@>>>G$ such that $F(D)=D$. Let $l$ be a 
prime number invertible in $\kk$; let $\bbq$ be an algebraic closure of the $l$-adic 
numbers. All group representations are assumed to be finite dimensional over $\bbq$. We say
"local system" instead of "$\bbq$-local system".

Let $\cb$ be the variety of Borel subgroups of $G^0$. Now $F:G@>>>G$ induces a morphism
$\cb@>>>\cb$ denoted again by $F$. We fix $B^*\in\cb$ and a maximal torus $T$ of $B^*$ such
that $F(B^*)=B^*$, $F(T)=T$. Let $U^*$ be the unipotent radical of $B^*$. Let $NB^*$ (resp.
$NT$) be the normalizer of $B^*$ (resp. $T$) in $G$. Let $\tT=NT\cap NB^*$, a closed 
$F$-stable subgroup of $G$ with identity component $T$. Let $\tT_D=\tT\cap D$. 

Let $\cn=NT\cap G^0$. Let $W=\cn/T$ be the Weyl group. Let $\uD:T@>\si>>T$, $\uD:W@>\si>>W$
be the automorphisms induced by $\Ad(d):\cn@>>>\cn$ where $d$ is any element of $\tT_D$. 
Now $F:\cn@>>>\cn$ 
induces an automorphism of $W$ denoted again by $F$. For $w\in W$ let $[w]$ be the inverse
image of $w$ under the obvious map $\cn@>>>W$ and let $\uw$ be the 
automorphism $\Ad(x):T@>>>T$ for any $x\in[w]$. For $w\in W$ let $\co_w$ be the $G^0$-orbit
in $\cb\T\cb$ ($G^0$ acting by simultaneous conjugation on both factors) that contains 
$(B^*,xB^*x\i)$ for some/any $x\in[w]$. Define the "length function" $l:W@>>>\NN$ by 
$l(w)=\dim\co_w-\dim\cb$. For any $y\in G^0$ we define $k(y)\in\cn$ by $y\in U^*k(y)U^*$. 
For $y\in G^0,\t\in\tT$ we have $k(\t y\t\i)=\t k(y)\t\i$ and $F(k(y))=k(F(y))$. For 
$x\in G^0$ we define $F_x:G@>>>G$ by $F_x(g)=xF(g)x\i$; this is the
Frobenius map for an $\FF_q$-rational structure on $G$. (Indeed if $y\in G^0$ is such that
$x=y\i F(y)$, then $\Ad(y):G@>\sim>>G$ carries $F_x$ to $F$.) If $w\in W$ satisfies 
$\uD(w)=w$ and $x\in[w]$ then $T,\tT$ are $F_x$-stable; thus $F_x$ is the Frobenius map for
an $\FF_q$-rational structure on $\tT$ whose group of rational points is $\tT^{F_x}$. Since
$\tT_D^{F_x}$ is the set of rational points of $\tT_D$ (a homogeneous $T$-space under left
translation) for the rational structure defined by $F_x:\tT_D@>>>\tT_D$, we have 
$\tT_D^{F_x}\ne\em$. 

Let $Z_\em=\{(B_0,g)\in\cb\T D;gB_0g\i=B_0\}$. Let $d\in\tT_D$. We set
$$\dZ_{\em,d}=\{(h_0U^*,g)\in(G^0/U^*)\T D;h_0\i gh_0d\i\in B^*\}.$$   
Define $a_\em:\dZ_{\em,d}@>>>Z_\em$ by $(h_0U^*,g)\m(h_0B^*h_0\i,g)$. Now $a_\em$ is a 
principal $T$-bundle where $T$ acts (freely) on $\dZ_{\em,d}$ by
$t_0:(h_0U^*,g)\m(h_0t_0\i,g)$. Define $p_\em:Z_\em@>>>D$ by $(B_0,g)\m g$. We define 
$b_\em:\dZ_{\em,d}@>>>T$ by $(h_0U^*,g)\m k(h_0\i gh_0d\i)$. Note that $b_\em$ commutes 
with the $T$-actions where $T$ acts on $T$ by 

(a) $t_0:t\m t_0t\uD(t_0\i)$.
\nl
Let $\cl$ be a local system of rank $1$ on $T$ such that

(i) $\cl^{\ot n}\cong\bbq$ for some $n\ge1$ invertible in $\kk$;

(ii) $\uD^*\cl\cong\cl$;
\nl
From (i),(ii) we see (using \cite{\CDG, 28.2(a)}) that $\cl$ is equivariant for the 
$T$-action (a) on $T$. Hence $b_\em^*\cl$ is a $T$-equivariant local system on 
$\dZ_{\em,d}$. Since $a_\em$ is a principal $T$-bundle there is a well defined local 
system $\tcl_\em$ on $Z_\em$ such that $a_\em^*\tcl_\em=b_\em^*\cl$. Note that the 
isomorphism class of $\tcl_\em$ is independent of the choice of $d$. Assume in addition 
that:

(iii) $\{w\in W;\uD(w)=w,\uw^*\cl\cong\cl\}=\{1\}$.
\nl
We show:

(b) {\it $p_{\em!}\tcl_\em$ is an irreducible intersection cohomology complex on $D$.}
\nl
We identify $Z_\em$ with the variety $X=\{(g,xB^*)\in G\T G^0/B^*;x\i gx\in NB^*\}$ (as in
\cite{\CDG, I, 5.4} with $P=B^*,L=T,S=\tT_D$) by $(g,xB^*)\lra(xB^*x\i,g)$. Then $\tcl_\em$
becomes the local system $\bar\ce$ on $X$ defined as in \cite{\CDG, I, 5.6} in terms of the
local system $\ce=j^*\cl$ on $\tT_D$ where $j:\tT_D@>>>T$ is $y\m d\i y$. (Note that $\ce$
is equivariant for the conjugation action of $T$ on $\tT_D$.) In our case we have 
$\bar\ce=IC(X,\bar\ce)$ since $X$ is smooth. Hence from \cite{\CDG, I, 5.7} we see that 
$p_{\em!}\bar\ce$ is an intersection cohomology complex on $D$ corresponding to a 
semisimple local system on an open dense subset of $D$ which, by the results in 
\cite{\CDG, II, 7.10}, is irreducible if and only if the following condition is satisfied: 
if $w\in W,x\in[w]$ satisfy $\Ad(x)(\tT_D)=\tT_D$ and $\Ad(x)^*\ce\cong\ce$, then $w=1$. 
This is clearly equivalent to condition (iii). This proves (b).

From (b) and the definitions we see that $p_{\em!}\tcl_\em[\dim D]$ is a character sheaf on
$D$ in the sense of \cite{\CDG, VI}. A character sheaf on $D$ of this form is said to be 
{\it generic}. We can state the following result.

\proclaim{Theorem 1.2} Let $A$ be a generic character sheaf on $D$ such that $F^*A\cong A$
where $F:D@>>>D$ is the restriction of $F:G@>>>G$. Let $\ps:F^*A@>>>A$ be an isomorphism.
Define $\c_\ps:D^F@>>>\bbq$ by $g\m\sum_{i\in\ZZ}(-1)^i\tr(\ps,\ch^i_g(A))$ where $\ch^i$ 
is the $i$-th cohomology sheaf and $\ch^i_g$ is its stalk at $g$. There exists a 
$G^F$-module $V$ and a scalar $\l\in\bbq^*$ such that $\c_\ps(g)=\l\tr(g,V)$ for all 
$g\in D^F$.
\endproclaim
The proof is given in \S3. We now make some preliminary observations. In the setup of
1.1 we have $A=p_{\em!}\tcl_\em[\dim D]$ where $\cl$ satisfies 1.1(i),(ii),(iii) and 
$F^*(p_{\em!}\tcl_\em)\cong p_{\em!}\tcl_\em$. Hence we have 
$p_{\em!}\wt{F^*\cl}_\em\cong p_{\em!}\tcl_\em$. By a computation in \cite{\CDG, IV, 21.18}
we deduce that there exists $w'\in W$ such that $\uD(w')=w'$, $\uw'{}^*F^*\cl\cong\cl$. 
Setting $w=F(w')$ we see that 

(a) $\uD(w)=w$, $F^*\uw^*\cl\cong\cl$.

\subhead 1.3\endsubhead
Let $\ww=(w_1,w_2,\do,w_r)$ be a sequence in $W$. Let $l_\ww=l(w_1)+l(w_2)+\do+l(w_r)$. Let
$$Z_\ww=\{(B_0,B_1,\do,B_r,g)\in\cb^{r+1}\T D;
gB_0g\i=B_r,(B_{i-1},B_i)\in\co_{w_i}(i\in[1,r])\}.$$
This agrees with the definition in 1.1 when $r=0$, that is $\ww=\em$. Let $d\in\tT_D$. We 
define $\dZ_{\ww,d}$ as in 1.1 when $r=0$ and by
$$\align\dZ_{\ww,d}&=\{(h_0U^*,h_1B^*,\do,h_{r-1}B^*,h_rU^*,g)\in\\&
(G^0/U^*)\T(G^0/B^*)\T\do\T(G^0/B^*)\T(G^0/U^*)\T D;\\&
k(h_{i-1}\i h_i)\in[w_i](i\in[1,r]),h_r\i gh_0d\i\in U^*\};\endalign$$   
when $r\ge1$. Define $a_\ww:\dZ_{\ww,d}@>>>Z_\ww$ as in 1.1 when $r=0$ and by
$$\align&(h_0U^*,h_1B^*,\do,h_{r-1}B^*,h_rU^*,g)\m \\&
(h_0B^*h_0\i,h_1B^*h_1\i,\do,h_{r-1}B^*h_{r-1},h_rB^*h_r\i,g),\endalign$$
when $r\ge1$. Note that $a_\ww$ is a principal $T$-bundle where $T$ acts (freely) on 
$\dZ_{\ww,d}$ as in 1.1 when $r=0$ and by
$$\align&t_0:(h_0U^*,h_1B^*,\do,h_{r-1}B^*,h_rU^*,g)\m\\&
(h_0t_0\i U^*,h_1B^*,\do,h_{r-1}B^*,h_rdt_0\i d\i U^*,g)\endalign$$
when $r\ge1$. Define $p_\ww:Z_\ww@>>>D$ by $(B_0,B_1,\do,B_r,g)\m g$. 

In the remainder of this subsection we assume that $w_1w_2\do w_r=1$; this holds 
automatically when $r=0$. We define $b_\ww:\dZ_{\ww,d}@>>>T$ as in 1.1 when $r=0$ and by
$$(h_0U^*,h_1B^*,\do,h_{r-1}B^*,h_rU^*,g)\m k(h_0\i h_1)k(h_1\i h_2)\do k(h_{r-1}\i h_r)$$
when $r\ge1$. Note that $b_\ww$ commutes with the $T$-actions where $T$ acts on $T$ as in 
1.1(a).

Let $\cl$ be a local system of rank $1$ on $T$ such that 1.1(i),(ii) hold. As in 1.1, $\cl$
is equivariant for the $T$-action 1.1(a) on $T$. Hence $b_\ww^*\cl$ is a $T$-equivariant 
local system on $\dZ_{\ww,d}$. Since $a_\ww$ is a principal $T$-bundle there is a well 
defined local system $\tcl_\ww$ on $Z_\ww$ such that $a_\ww^*\tcl_\ww=b_\ww^*\cl$. 

\proclaim{Lemma 1.4} Assume that $w_1w_2\do w_r=1$ and that $\cl$ (as in 1.3) satisfies 

(i) $\cha^*\cl\not\cong\bbq$ for any coroot $\cha:\kk^*@>>>T$.
\nl
Then $p_{\ww!}\tcl_\ww[l_\ww](l_\ww/2)\cong p_{\em!}\tcl_\em$. (Note that $l_\ww$ is even.)
\endproclaim
Assume first that for some $i\in[1,r]$ we have $w_i=w'_iw''_i$ where $w'_i,w''_i$ in $W$
satisfy $l(w'_iw''_i)=l(w'_i)+l(w''_i)$. Let 
$$\ww'=(w_1,w_2,\do,w_{i-1},w'_i,w''_i,w_{i+1},\do,w_n).$$
The map $(B_0,B_1,\do,B_{r+1},g)\m(B_0,B_1,B_{i-1},B_{i+1},\do,B_{r+1},g)$ defines an 
isomorphism $Z_{\ww'}@>>>Z_\ww$ compatible with the maps $p_{\ww'},p_\ww$ and with the
local systems $\tcl_{\ww'},\tcl_\ww$. Since $l_{\ww'}=l_\ww$ we have

(a) $p_{\ww!}\tcl_\ww[l_\ww](l_\ww/2)\cong p_{\ww'!}\tcl_{\ww'}[l_{\ww'}](l_{\ww'}/2)$.
\nl
Using (a) repeatedly we can assume that $l(w_i)=1$ for all $i\in[1,r]$. We will prove the 
result in this case by induction on $r$. Note that $r$ is even. When $r=0$ the result is 
obvious. We now assume that $r\ge2$. Since $w_1w_2\do w_r=1$, we can find $j\in[1,r-1]$ 
such that $l(w_1w_2\do w_j)=j$, $l(w_1w_2\do w_{j+1})=j-1$. We can find a sequence 
$\ww'=(w'_1,w'_2,\do,w'_r)$ in $W$ such that $l(w'_i)=1$ for all $i\in[1,r]$, 
$w'_1w'_2\do w'_j=w_1w_2\do w_j$, $w'_j=w'_{j+1}$, $w'_i=w_i$ for $i\in[j+1,r]$. Let 

$\uu=(w_1w_2\do w_j, w_{j+1},\do,w_r)=(w'_1w'_2\do w'_j, w'_{j+1},\do,w'_r)$.
\nl
Using (a) repeatedly we see that 
$$p_{\ww!}\tcl_\ww[l_\ww](l_\ww/2)\cong p_{\uu!}\tcl_\uu[l_\uu](l_\uu/2)
\cong p_{\ww'!}\tcl_{\ww'}[l_{\ww'}](l_{\ww'}/2).$$
Replacing $\ww$ by $\ww'$ we see that we may assume in addition that $w_j=w_{j+1}$ for some
$j\in[1,r-1]$. We have a partition $Z_\ww=Z'_\ww\cup Z''_\ww$ where $Z'_\ww$ (resp. 
$Z''_\ww$) is defined by the condition $B_{j-1}=B_{j+1}$ (resp. $B_{j-1}\ne B_{j+1}$). Let 
$\ww'=(w_1,w_,\do,w_{j-1},w_{j+2},\do,w_r)$, $\ww''=(w_1,w_,\do,w_{j-1},w_{j+1},\do,w_r)$.
Define $c:Z'_\ww@>>>Z_{\ww'}$ by 
$$(B_0,B_1,\do,B_r,g)\m(B_0,B_1,\do,B_{j-1},B_{j+2},\do,B_r,g).$$
This is an affine line bundle and $\tcl_\ww|_{Z'_\ww}=c^*\tcl_{\ww'}$. Let $p'_\ww$ be the
restriction of $p_\ww$ to $Z'_\ww$. We have $p'_\ww=p_{\ww'}c$. Since the induction 
hypothesis applies to $\ww'$ we have 
$$\align&p'_{\ww!}(\tcl_\ww|_{Z'_\ww})[l_\ww](l_\ww/2)
=p_{\ww'!}c_!c^*\tcl_{\ww'}[l_\ww](l_\ww/2)\\&=
p_{\ww'!}\tcl_{\ww'}[-2](-1)[l_\ww](l_\ww/2)=p_{\ww'!}\tcl_{\ww'}[l_{\ww'}](l_{\ww'}/2)
=p_{\em!}\tcl_\em.\tag b\endalign$$
Define $e:Z''_\ww@>>>Z_{\ww''}$ by 
$$(B_0,B_1,\do,B_r,g)\m(B_0,B_1,\do,B_{j-1},B_{j+1},\do,B_r,g).$$
Let $p''_\ww$ be the restriction of $p_\ww$ to $Z''_\ww$. We have $p''_\ww=p_{\ww''}e$. We
show that $p''_{\ww!}(\tcl_\ww|_{Z''_\ww})=0$. It is enough to show that
$$p_{\ww''!}e_!(\tcl_\ww|_{Z''_\ww})=0.\tag c$$
Hence it is enough to show that $e_!(\tcl_\ww|_{Z''_\ww})=0$. It is also enough to show 
that, if $E$ is a fibre of $e$, then $H^i_c(E,\tcl_\ww|_E)=0$ for any $i$. As in the proof
of \cite{\CDG, VI, 28.10} we may identify $E=\kk^*$ in such a way that $\tcl_\ww|_E$ 
becomes $\cha^*(\cl)$ for some coroot $\cha:\kk^*@>>>T$. We then use that 
$H^i_c(\kk^*,\cha^*\cl)=0$ which follows from $\cha^*\cl\not\cong\bbq$.

Using (c) and the exact triangle 
$$(p_{\ww''!}e_!(\tcl_\ww|_{Z''_\ww}),p_{\ww!}\tcl_\ww,p'_{\ww!}(\tcl_\ww|_{Z'_\ww}))$$
we see that
$$p_{\ww!}\tcl_\ww[l_\ww](l_\ww/2)=p'_{\ww!}(\tcl_\ww|_{Z'_\ww})[l_\ww])(l_\ww/2)
=p_{\em!}\tcl_\em$$
(the last equality follows from (b)). The lemma is proved.

\proclaim{Lemma 1.5} Assume that $\cl$ (as in 1.3) satisfies 1.1(iii). Then $\cl$ satisfies
1.4(i).
\endproclaim
Let $R_\cl$ be the set of roots $\a:T@>>>\kk^*$ such that the corresponding coroot $\cha$ 
satisfies $\cha^*\cl\cong\bbq$. Let $W_\cl$ be the subgroup of $W$ generated by the
reflections with respect to the various $\a\in R_\cl$. Since $\uD^*\cl\cong\cl$ we have 
$\uD(W_\cl)=W_\cl$. Assume that 1.4(i) does not 
hold. Then $R_\cl\ne\em$ and $W_\cl\ne\{1\}$. By \cite{\DL, 5.17} the fixed point set of 
$\uD:W_\cl@>>>W_\cl$ is $\ne\{1\}$. Let $w\in W_\cl-\{1\}$ be such that $\uD(d)w=w$. Since
$w\in W_\cl$ we have $\uw^*\cl\cong\cl$ (see \cite{\CDG, VI, 28.3(b)}). Thus 1.1(iii) does 
not hold. The lemma is proved.

\head 2. Constructing representations of $G^F$\endhead
\subhead 2.1\endsubhead
In this section we construct some representations of $G^F$ using the method of \cite{\DL}.
See \cite{\MA},\cite{\DM} for other results in this direction.

Let $\cl$ be a local system of rank $1$ on $T$ such that 1.1(i) holds. For any $t\in T$ let
$\cl_t$ be the stalk of $\cl$ at $t$. Assume that we are given $w\in W$ and $x\in[w]$ such
that

(i) $F_x^*\cl\cong\cl$;
\nl
($F_x:T@>>>T$ as in 1.1). Let $\ph:F_x^*\cl@>>>\cl$ be the unique isomorphism of local 
systems on $T$ which induces the identity map on $\cl_1$. For $t\in T$, $\ph$ induces an 
isomorphism $\cl_{F_x(t)}@>\si>>\cl_t$. When $t\in T^{F_x}$ this is an automorphism of the 
$1$-dimensional vector space $\cl_t$ given by multiplication by $\th(t)\in\bbq^*$. It is 
well known that $t\m\th(t)$ is a group homomorphism $T^{F_x}@>>>\bbq^*$.

Following \cite{\DL} we define
$$Y=\{hU^*\in G^0/U^*;h\i F(h)\in U^*xU^*\}.$$
For $(g,t)\in G^{0F}\T T^{F_x}$ we define $e_{g,t}:Y@>>>Y$ by $hU^*\m ght\i U^*$. Note that
$(g,t)\m e_{g,t}$ is an action of $G^{0F}\T T^{F_x}$ on $Y$. Hence $G^{0F}\T T^{F_x}$ acts
on $H^i_c(Y):=H^i_c(Y,\bbq)$ by $(g,\t)\m e_{g\i,\t\i}^*$. We set
$$H^i_c(Y)_\th=\{\x\in H^i_c(Y);e_{1,t\i}^*\x=\th(t)\i\x\text{ for all }t\in T^{F_x}\};$$ 
this is a $G^{0F}\T T^{F_x}$-stable subspace of $H^i_c(Y)$.

For $g\in G^{0F}$ we define $\e_g:H^i_c(Y)_\th@>>>H^i_c(Y)_\th$ by $\e_g(\x)=e_{g\i,1}^*$. 
This makes $H^i_c(Y)_\th$ into a $G^{0F}$-module.

We can find an integer $r\ge1$ such that 
$$F^r(x)=x,\qua xF(x)\do F^{r-1}(x)=1.$$ 
Indeed we first find an integer $r_1\ge1$ such that $F^{r_1}(x)=x$ and then we find an 
integer $r_2\ge1$ such that $(xF(x)\do F^{r_1-1}(x))^{r_2}=1$. Then $r=r_1r_2$ has the 
required properties. Then $hU^*\m F^r(h)U^*$ is a well defined map $Y@>>>Y$ denoted again 
by $F^r$. Also,
$$F^r=F_x^r:G@>>>G.$$
(We have 
$F_x^r(g)=(xF(x)\do F^{r-1}(x))F^r(g)(xF(x)\do F^{r-1}(x))\i=F^r(g)$.)
Hence $F^r$ acts trivially on $T^{F_x}$. We see that $F^r:Y@>>>Y$ commutes with 
$e_{g,t}:Y@>>>Y$ for any $(g,t)\in G^{0F}\T T^{F_x}$. Hence $(F^r)^*:H^i_c(Y)@>>>H^i_c(Y)$
leaves stable the subspace $H^i_c(Y)_\th$. Note that:

for any $i$, all eigenvalues of $(F^r)^*:H^i_c(Y)@>>>H^i_c(Y)$ are of the form root of $1$
times $q^{nr/2}$ where $n\in\ZZ$.
\nl
(See \cite{\GF, 6.1(e)} and the references there.) 

Replacing $r$ by an integer multiple we may therefore assume that $r$ satisfies in addition
the following condition:

(a) for any $i$, all eigenvalues of $(F^r)^*:H^i_c(Y)@>>>H^i_c(Y)$ are of the form 
$q^{nr/2}$ where $n\in\ZZ$.

\subhead 2.2\endsubhead
We preserve the setup of 2.1 and assume in addition that $\cl$ satisfies 1.4(i). 
Let $i_0=2\dim U^*-l(w)$. Note that

(a) {\it $H^i_c(Y)_\th=0$ for $i\ne i_0$; if $i=i_0$ then all 
eigenvalues of $(F^r)^*:H^i_c(Y)_\th@>>>H^i_c(Y)_\th$ are of the form $q^{ir/2}$.}
\nl
For the first statement in (a) see \cite{\DL, 9.9} and the remarks in the proof of
\cite{\GF, 8.15}. The second statement in (a) is deduced from 2.1(a) as in the proof of
\cite{\GF, 6.6(c)}.

\subhead 2.3\endsubhead
We preserve the setup of 2.1 and assume in addition that $\cl$ satisfies 1.1(ii) and that 
$w\in W$ satisfies $\uD(w)=w$. From the definitions we see that $\uD:T@>>>T$ commutes with
$F_x:T@>>>T$ hence $\uD$ restricts to an automorphism of $T^{F_x}$ and that

(a) $\th(\uD(t))=\th(t)$ for any $t\in T^{F_x}$.
\nl
We show:

(b) {\it there exists a homomorphism $\ti\th:\tT^{F_x}@>>>\bbq^*$ such that 
$\ti\th|_{T^{F_x}}=\th$.}
\nl
Let $d\in\tT^{F_x}_D$. Let $n=|G/G^0|=|\tT^{F_x}/T^{F_x}|$. Then $t_0:=d^n\in T^{F_x}$. Let
$c\in\bbq^*$ be such that $c^n=\th(t_0)$. For any $t\in T^{F_x}$ and $j\in\ZZ$ we set 
$\ti\th(d^jt)=c^j\th(t)$. This is well defined: if $d^jt=d^{j'}t'$ with $j,j'\in\ZZ$, 
$t,t'\in T^{F_x}$ then $j'=j+nj_0$, $j_0\in\ZZ$ and $t'=t_0^{j_0}t$ so that 
$\th(t')=c^{nj_0}\th(t)$ and $c^j\th(t)=c^{j'}\th(t')$. We show that if $j,j'\in\ZZ$, 
$t,t'\in T^{F_x}$ then $\ti\th(d^jtd^{j'}t')=\ti\th(d^jt)\ti\th(d^{j'}t')$ that is
$c^{j+j'}\th(\uD^{-j'}(t)t')=c^j\th(t)c^{j'}\th(t')$; this follows from (a). This proves 
(b).

Let $\G=\{(g,\t)\in G^F\T\tT^{F_x};g\t\i\in G^0\}$, a subgroup of $G^F\T\tT^{F_x}$. For
$(g,\t)\in\G$ we define $e_{g,\t}:Y@>>>Y$ by $hU^*\m gh\t\i U^*$. To see that this is well
defined we assume that $h\in G^0$ satisfies $h\i F(h)\in U^*xU^*$ and $(g,\t)\in\G$; we 
compute
$$\align&(gh\t\i)\i F(gh\t\i)=\t h\i g\i gF(h)F(\t\i)
\\&=\t h\i F(h)F(\t\i)\in\t U^*xU^*F(\t\i)=U^*\t xF(\t\i)U^*=U^*xU^*,\endalign$$
since $\t xF(\t\i)=x$ (that is $F_x(\t)=\t$). Note that $(g,\t)\m e_{g,\t}$ is an action of
$\G$ on $Y$ (extending the action of $G^{0F}\T T^{F_x}$). Hence $\G$ acts on $H^i_c(Y)$ by
$(g,\t)\m e_{g\i,\t\i}^*$. Note that $H^i_c(Y)_\th$ is a $\G$-stable subspace of 
$H^i_c(Y)$. This follows from the identity 
$$e_{g\i,\t\i}e_{1,t\i}=e_{1,\t\i t\i\t}e_{g\i,\t\i}$$
for $g\in G^F$, $\t\in\tT^{F_x}$, $t\in T^{F_x}$ together with the identity 
$\th(t)=\th(\t\i t\t)$ which is a consequence of (a).

For $g\in G^F$ we define $\e_g:H^i_c(Y)_\th@>>>H^i_c(Y)_\th$ by
$$\e_g(\x)=\ti\th(\t)e_{g\i,\t\i}^*\x$$ 
for any $\x\in H^i_c(Y)_\th$  and any $\t\in\tT^{F_x}$ such that $g\t\i\in G^0$. Assume 
that $\t'\in\tT^{F_x}$ is another element such that $g\t'{}\i\in G^0$. Then $\t'=\t t$ with
$t\in T^{F_x}$ and 
$$\ti\th(\t')e_{g\i,\t'{}\i}^*\x=\ti\th(\t)\th(t)e_{g\i,\t{}\i}^* e_{1,t\i}^*\x=
\ti\th(\t)e_{g\i,\t{}\i}^*\x$$
so that $\e_g$ is well defined. For $g,g'$ in $G^F$ we choose $\t,\t'$ in $\tT^{F_x}$ such
that $g\t\i\in G^0,g'\t'{}\i\in G^0$; we have
$$\e_g\e_{g'}\x=\ti\th(\t')\ti\th(\t)e_{g\i,\t\i}^*e_{g'{}\i,\t'{}\i}^*\x=
\ti\th(\t\t')e_{(gg')\i,(\t\t')\i}^*\x=\e_{gg'}\x.$$
We see that

{\it $g\m\e_g$ defines a $G^F$-module structure on $H^i_c(Y)_\th$ extending the
$G^{0F}$-module structure in 2.1.}
\nl
(Note that this extension depends on the choice of $\ti\th$.) We show:

(c) {\it If $(g,\t)\in\G$ then $F^re_{g,\t}:Y@>>>Y$ is the Frobenius map of an 
$\FF_q$-rational structure on $Y$.}
\nl
Since $e_{g,t}$ is a part of a $\G$-action, it has finite order. Since $F^r=F_x^r:G@>>>G$
(see 2.1), we see that $F^r:Y@>>>Y$ commutes with $e_{g,\t}:Y@>>>Y$. Hence (c) holds.

\subhead 2.4\endsubhead
We preserve the setup of 2.3 and assume in addition that $\cl$ satisfies 1.3(i). Let 
$i_0=2\dim U^*-l(w)$. Using 2.2(a), 2.3(c) and Grothendieck's trace formula we see that for
$(g,d)\in\G$ we have 
$$\align&(-1)^{l(w)}\ti\th(d)q^{i_0r/2}\tr(\e_g,H^{i_0}_c(Y)_\th)\\&=
\ti\th(d)\sum_i(-1)^i\tr((F^r)^*\e_g,H^i_c(Y)_\th)=\sum_i(-1)^i
\tr((F^r)^*e_{g\i,d\i}^*,H^i_c(Y)_\th)\\&=\sum_i(-1)^i|T^{F_x}|\i\sum_{t\in T^{F_x}}
\tr((F^r)^*e_{g\i,d\i}^*e_{1,t\i}^*,H^i_c(Y))\th(t)\\&
=|T^{F_x}|\i\sum_{t\in T^{F_x}}\sum_i(-1)^i\tr((F^r)^*e_{g\i,(dt)\i}^*,H^i_c(Y))\th(t)\\&
=|T^{F_x}|\i\sum_{t\in T^{F_x}}|Y^{F^re_{g\i,(dt)\i}}|\th(t)\\&=|T^{F_x}|\i\sum_{t\in T^
{F_x}}|\{hU^*\in(G^0/U^*);h\i F(h)\in U^*xU^*,h\i g\i F^r(h)dt\in U^*\}|\th(t).\endalign$$

\head 3. Proof of Theorem 1.2\endhead
\subhead 3.1\endsubhead
Let $A,\ps,\c_\ps$ be as in 1.2. Let $\cl,w$ be as in the end of 1.2. Let $x\in[w]$. From 
1.2(a) we see that 2.1(i) holds. Let $r\ge1$ be as in 2.1. Let 
$$\ww=(w,F(w),\do,F^{r-1}(w)).$$ 
By the choice of $r$ we have $wF(w)\do F^{r-1}(w)=1$. Define a morphism 
$\tF:Z_\ww@>>>Z_\ww$ by 
$$\tF(B_0,B_1,\do,B_r,g)=(F(g\i B_{r-1}g),F(B_0),F(B_1),\do,F(B_{r-1}),F(g)).$$
We show:

(a) {\it Let $g\in D^F$ and let $\tF_g:p_\ww\i(g)@>>>p_\ww\i(g)$ be the restriction of 
$\tF:Z_\ww@>>>Z_\ww$. Then $\tF_g$ is the Frobenius map of an $\FF_q$-rational structure on
$p_\ww\i(g)$.}
\nl
It is enough to note that the map $\cb^{r+1}@>>>\cb^{r+1}$ given by
$$(B_0,B_1,\do,B_r)\m(F(g\i B_{r-1}g),F(B_0),F(B_1),\do,F(B_{r-1}))$$ 
is the composition of the map 
$$F':(B_0,B_1,\do,B_r)\m(F(B_0),F(B_1),\do,F(B_r))$$ 
(the Frobenius map of an $\FF_q$-rational structure on $\cb^{r+1}$) with the automorphism 
$$(B_0,B_1,\do,B_r)\m(g\i B_{r-1}g,B_0,B_1,\do,B_{r-1})$$ 
of $\cb^{r+1}$ which commutes with $F'$ and has finite order (since $g$ has finite order 
in $G$).

Let $d\in\tT^{F_x}_D$. Define a morphism $\tF':\dZ_{\ww,d}@>>>\dZ_{\ww,d}$ by 
$$\tF'(h_0U^*,h_1B^*,\do,h_{r-1}B^*,h_rU^*,g)
=(h'_0U^*,h'_1B^*,\do,h'_{r-1}B^*,h'_rU^*,F(g))$$
where
$$h'_0=F(g\i h_{r-1}k(h_{r-1}\i h_r))x\i d,\qua h'_r=F(h_{r-1}k(h_{r-1}\i h_r)x\i,$$
$$h'_i=F(h_{i-1})\text{ for }i\in[1,r-1].$$
This is well defined since 
$$(F(h_{r-1}k(h_{r-1}\i h_r)x\i)\i F(g)F(g\i h_{r-1}k(h_{r-1}\i h_r))x\i)dd\i=1.$$

We show that the $T$-action on $\dZ_{\ww,d}$ (see 1.3) satisfies 
$\tF'(t_0\tx)=F_x(t_0)\tF'(\tx)$ for $t_0\in T,\tx\in\dZ_{\ww,d}$. Let $(h_i)$ be as 
above. We must show:
$$F(g\i h_{r-1}k(h_{r-1}\i h_rdt_0\i d\i))x\i d=
F(g\i h_{r-1}k(h_{r-1}\i h_r))x\i dxF(t_0\i)x\i,$$
$$F(h_{r-1}k(h_{r-1}\i h_rdt_0\i d\i)x\i=F(h_{r-1}k(h_{r-1}\i h_r)x\i dxF(t_0)\i x\i d\i,$$
which follow from $F(d)=x\i dx$. Note that 

(b) $a_\ww\tF'=\tF a_\ww:\dZ_{\ww,d}@>>>Z_\ww$. 
\nl
We show:

(c) {\it $|a_\ww\i(y)^{\tF'}|=|T^{F_x}|$ for any $y\in Z_\ww^{\tF}$.}
\nl
Since $a_\ww\i(y)$ is a homogeneous $T$-space this follows from Lang's theorem applied to
$(T,F_x)$.

We have 

(d) $p_\ww \tF=F p_\ww:Z_\ww@>>>D$.

\subhead 3.2\endsubhead
We show:

(a) $b_\ww\tF'=F_xb_\ww:\dZ_{\ww,d}@>>>T$.
\nl
Let $(h_0,h_1,\do,h_r,g)\in(G^0)^{r+1}\T D$ be such that 
$$(h_0U^*,h_1B^*,\do,h_{r-1}B^*,h_rU^*,g)\in\dZ_{\ww,d}.$$ 
Let $(h'_1,h'_2,\do,h'_r)$ be as in 3.1. We set
$$\mu=k(h_0\i h_1)k(h_1\i h_2)\do k(h_{r-1}\i h_r)\in T,$$
$$\mu'=k(h_0\i h_1)k(h_1\i h_2)\do k(h_{r-2}\i h_{r-1})\in B^*F^{r-1}(x)\i B^*$$
$$\ti\mu=k(h'_0{}\i h'_1)k(h'_1{}\i h'_2)\do k(h'_{r-1}{}\i h'_r)\in T$$
so that $\mu=\mu'k(h_{r-1}\i h_r)$ and 
$$\align&\ti\mu=k(d\i xF(k(h_{r-1}\i h_r)\i h_{r-1}\i gh_0))\\&
\T k(F(h_0\i h_1))\do k(F(h_{r-3}\i h_{r-2}))k(F(h_{r-2}\i h_{r-1}k(h_{r-1}\i h_r))x\i)\\&
=d\i xF(k(h_{r-1}\i h_r)\i)F(d)k(F(d\i)F(h_{r-1}\i gh_0))F(\mu')F(k(h_{r-1}\i h_r))x\i\\&
=d\i xF(d)F(\mu)x\i=xF(\mu)x\i=F_x(\mu),\endalign$$
as required. 

\subhead 3.3\endsubhead
Let $\ph:F_x^*\cl@>\si>>\cl$, $\th:T^{F_x}@>>>\bbq^*$ be as in 2.1. We shall denote by ? 
the various isomorphisms induced by $\ph$ such as:

(a) $\tF'{}^*b_\ww^*\cl=b_\ww^*F_x^*\cl@>\si>>b_\ww^*\cl$ (see 3.2(a)),

(b) $\tF'{}^*a_\ww^*\tcl_\ww@>\si>>a_\ww^*\tcl_\ww$ (coming from (a)),

(c) $a_\ww^*\tF^*\tcl_\ww@>\si>>a_\ww^*\tcl_\ww$ (see (b) and 3.1(b)),

(d) $\tF^*\tcl_\ww@>\si>>\tcl_\ww$ (coming from (c)),

(e) $p_{\ww!}\tF^*\tcl_\ww@>\si>>p_{\ww!}\tcl_\ww$ (coming from (d)),

(f) $F^*p_{\ww!}\tcl_\ww@>\si>>p_{\ww!}\tcl_\ww$ (coming from (e) and 3.1(d)).

(g) $F^*(p_{\ww!}\tcl_\ww[l_\ww])@>\si>>p_{\ww!}\tcl_\ww[l_\ww]$ (coming from (f)).

\subhead 3.4\endsubhead
For any $g\in D^F$ we compute 
$$\align&\sum_i(-1)^i\tr(?,\ch^i_g(p_{\ww!}\tcl_\ww))=\sum_i(-1)^i\tr(?,H^i_c(p_\ww\i(g),
\tcl_\ww))\\&=\sum_{y\in p_\ww\i(g);\tF(y)=y}\tr(?,(\tcl_\ww)_y)\endalign$$
where $\ch^i$ is the $i$-th cohomology sheaf. (The last two sums are equal by the 
Grothendieck trace formula applied in the context of 3.1(a).) Using 3.1(c) we see that the
last sum equals
$$\align&|T^{F_x}|\i\sum_{\ty\in a_\ww\i(p_\ww\i(g))^{\tF'}}\tr(?,(a_\ww^*\tcl_\ww)_{\ty})
=|T^{F_x}|\i\sum_{\ty\in a_\ww\i(p_\ww\i(g))^{\tF'}}\tr(?,(b_\ww^*\cl_\ww)_{\ty})\\&
=|T^{F_x}|\i\sum_{\ty\in a_\ww\i(p_\ww\i(g))^{\tF'}}\tr(?,(\cl_\ww)_{b_\ww(\ty)}).\endalign
$$
Now $a_\ww\i(p_\ww\i(g))^{\tF'}$ can be identified with the set of all
$$(h_0U^*,h_1B^*,\do,h_{r-1}B^*,h_rU^*)\in(G^0/U^*)\T(G^0/B^*)\T\do\T(G^0/B^*)\T(G^0/U^*)$$
such that

(a) $k(h_{i-1}\i h_i)\in F^{i-1}(x)T$ for $i\in[1,r]$,

(b) $h_r\i gh_0d\i\in U^*$,

(c) $h_0U^*=F(g\i h_{r-1}k(h_{r-1}\i h_r))x\i dU^*$, 

(d) $h_iB^*=F(h_{i-1})B^*$ for $i\in[1,r-1]$.
\nl
(We then have automatically $h_rU^*=F(h_{r-1}k(h_{r-1}\i h_r)x\i U^*$.) If $h_0U^*$ is 
given, then (d) determines successively $h_2B^*,\do h_{r-1}B^*$ in a unique way and (b) 
determines $h_rU^*$ in a unique way. We see that the equations (a)-(d) are equivalent to 
the following equations for $h_0U^*$:
$$h_0\i F(h_0)\in B^*xB^*,\qua F^{r-1}(h_0)\i gh_0d\i\in B^*F^{r-1}(x)B^*,$$
$$F^r(h_0)\i gh_0d\i U^*=k(F^r(h_0)\i gF(h_0)F(d\i))x\i U^*$$
(if $r\ge2$) and
$$h_0\i gh_0d\i\in B^*xB^*,\qua F(h_0)\i gh_0d\i U^*=k(F(h_0)\i gF(h_0)F(d\i))x\i U^*$$
(if $r=1$). In both cases these equations are equivalent to
$$h_0\i F(h_0)\in U^*txF(t)\i U^*,\qua F^r(h_0)\i gh_0d\i\in F^r(t)U^*\tag e$$
for some $t\in T$. We then have $F^{r-1}(h_0)\i gh_0d\i\in U^*F^{r-1}(t)F^{r-1}(x)U^*$. For
$h_0U^*,t$ as in (e) we compute
$$\align&k(h_0\i F(h_0))k(F(h_0)\i F^2(h_0))\do k(F^{r-2}(h_0)\i F^{r-1}(h_0)) 
                                                  k(F^{r-1}(h_0)\i gh_0d\i)\\&
=(txF(t)\i)(F(t)F(x)F^2(t\i))\do(F^{r-2}(t)F^{r-2}(x)F^{r-1}(t)\i)(F^{r-1}(t)F^{r-1}(x))
\\&=txF(x)\do F^{r-1}(x)=t.\endalign$$
By 3.2(a) the result of the last computation is necessarily in $T^{F_x}$. Thus $F_x(t)=t$. Hence $F^r(t)=t$ and the equations (e) become
$$h_0\i F(h_0)\in U^*xU^*,\qua F^r(h_0)\i gh_0d\i\in T^{F_x}U^*.\tag f$$
We see that  
$$\sum_i(-1)^i\tr(?,\ch^i_g(p_{\ww!}\tcl_\ww))=|T^{F_x}|\i\sum_{t\in T^{F_x}}a_t
=|T^{F_x}|\i\sum_{t'\in T^{F_x}}a'_{t'}$$
where
$$a_t=|\{hU^*\in(G^0/U^*);h\i F(h)\in U^*xU^*,dh\i g\i F^r(h)t\in U^*\}|\th(t),$$
$$a'_{t'}=|\{hU^*\in(G^0/U^*);h\i F(h)\in U^*xU^*,h\i g\i F^r(h)dt'\in U^*\}|\th(dt'd\i).$$
Comparing with the last formula in 2.4 and using $\th(dt'd\i)=\th(t')$ for $t'\in T^{F_x}$
we obtain (with $i_0$ as in 2.4):
$$\sum_i(-1)^i\tr(?,\ch^i_g(p_{\ww!}\tcl_\ww))=
(-1)^{l(w)}\ti\th(d)q^{i_0r/2}\tr(\e_g,H^{i_0}_c(Y)_\th).$$
Let us choose an isomorphism $p_{\ww!}\tcl_\ww[l_\ww]\cong p_{\em!}\tcl_\em$.
(This exists by 1.4; note that 1.4(i) holds by 1.5.) Via this isomorphism, the isomorphism
3.3(g) corresponds to an isomorphism $F^*(p_{\em!}\tcl_\em)@>>>p_{\em!}\tcl_\em$ that is to
an isomorphism $\ps':F^*A@>\sim>>A$ so that
$$\sum_i(-1)^i\tr(?,\ch^i_g(p_{\ww!}\tcl_\ww))=\sum_i(-1)^i\tr(\ps',\ch^i_g(A))$$
for any $g\in D^F$. (We use that $l_\ww$ is even.) Since $A$ is irreducible, we must have 
$\ps=\l'\ps'$ for some $\l'\in\bbq^*$. It follows that 
$$\sum_{i\in\ZZ}(-1)^i\tr(\ps,\ch^i_g(A))
=\l'(-1)^{l(w)}\ti\th(d)q^{i_0r/2}\tr(\e_g,H^{i_0}_c(Y)_\th)$$
for any $g\in D^F$. Thus Theorem 1.2 holds with $V$ being the $G^F$-module 
$H^{i_0}_c(Y)_\th$, which is irreducible (even as a $G^{0F}$-module) if $G^0$ has connected
centre, but is not necessarily irreducible in general.

\enddocument